\documentstyle[12pt,amssymb,amstex]{amsart}
\numberwithin{equation}{section} \textwidth 140mm \textheight 220mm
\def\cb{{\cal B}}

\def\cf{{\cal F}}

\def\ch{{\cal H}}

\def\ck{{\cal K}}

\def\co{{\cal O}}

\def\car{{\cal R}}
\def\cs{{\cal S}}

\def\ga{{\frak A}}

\def\gg{{\frak G}}

\def\gar{{\frak R}}

\def\bc{{\mathbb C}}

\def\bbf{{\mathbb F}}

\def\bn{{\mathbb N}}
\def\bp{{\mathbb P}}
\def\bq{{\mathbb Q}}
\def\br{{\mathbb R}}
\def\bt{{\mathbb T}}

\def\bz{{\mathbb Z}}

\def\a{\alpha}

 \def\G{\Gamma}
\def\d{\delta}

\def\l{\lambda} 

\def\m{\mu}

\def\n{\nu}
\def\r{\rho}
\def\s{\sigma} 

\def\f{\varphi} 
\def\th{\theta} \def\Th{\Theta}
\def\om{\omega} \def\Om{\Omega}

\newtheorem{thm}{Theorem}[section]
\newtheorem{lem}[thm]{Lemma}

\newtheorem{prop}[thm]{Proposition}
\newtheorem{defin}[thm]{Definition}

\def\di{\mathop{\rm d}}

\newcommand{\nn}{\nonumber}
\def\idd{{\bf 1}\!\!{\rm I}}

\begin{document}

\title[ergodic properties]
{ergodic properties of bogoliubov automorphisms in free probability}
\author{Francesco Fidaleo}
\address{Francesco Fidaleo \\Department of Comput. \&
Theor. Sci., Faculty of Science, IIUM, P.O. Box, 141, 25710,
Kuantan, Pahang, Malaysia. Permanent address: Dipartimento di
Matematica, II Universit\`{a} di Roma ``Tor Vergata'', Via della
Ricerca Scientifica, 00133 Roma, Italia} \email{{\tt
fidaleo@@mat.uniroma2.it}}
\author{Farrukh Mukhamedov}
\address{Farrukh Mukhamedov\\
Department of Comput. \& Theor. Sci., Faculty of Science, IIUM, P.O.
Box, 141, 25710, Kuantan, Pahang, Malaysia} \email{{\tt
far75m@@yandex.ru}, {\tt farrukh\_m@@iiu.edu.my}}

\begin{abstract}
We show that some $C^*$--dynamical systems obtained by "quantizing"
classical ones on the free Fock space, enjoy very strong ergodic
properties. Namely, if the classical dynamical system $(X, T, \m)$
is ergodic but not weakly mixing, then the resulting quantized
system $(\gg,\a)$ is uniquely ergodic (w.r.t the fixed point
algebra) but not uniquely weak mixing. The same happens if we
quantize a classical system $(X, T, \m)$ which is weakly mixing but
not mixing. In this case, the quantized system is uniquely weak
mixing but not uniquely mixing. Finally, a quantized system arising
from a classical mixing dynamical system, will be uniquely mixing.
In such a way, it is possible to exhibit uniquely weak mixing and
uniquely mixing $C^*$--dynamical systems whose GNS representation
associated to the unique invariant state generates a von Neuman
factor of one of the following types: $I_{\infty}$, $II_{1}$,
$III_{\l}$ where $\l\in(0,1]$. The results listed above are extended to the $q$--commutation relations, provided $|q|<\sqrt2-1$.
\end{abstract}

\maketitle

\section{introduction}

The study of quantum dynamical systems has been an impetuos growth
in the last years, in view of natural applications to various field
of mathematics and physics. It is then of interest to understand
among the various ergodic properties, which ones survive and are
meaningful, by passing from the classical to the quantum case. Due
to noncommutativity, in the latter situation the matter is much more
complicated than in the former. The reader is referred e.g. to
\cite{AH, J, NSZ} for further details relative to some differences
between classical and quantum situations. Therefore, it is then
natural to study of the possible generalizations to quantum case of
the various ergodic properties known for classical dynamical
systems.

By coming back to the classical case, one of the strong ergodic
properties of the dynamical system $(\Om,T)$ consisting of a compact
metric space $\Om$ and a homeomorphism $T$, is the unique ergodicity
which means that there exists a unique invariant Borel measure $\m$
for $T$. It is known (cf. \cite{KSF}) that the last property is
equivalent to the uniform convergence in $C(\Om)$ of the ergodic
averages ${\displaystyle\frac1n\sum_{k=0}^{n-1}f\circ T^{k}}$ to
$\int f\di\m$ for any $f\in C(\Om)$. A pivotal example of classical
uniquely ergodic dynamical system is given by an irrational rotation
on the unit circle (see \cite{Da} for further examples). In the
quantum setting, the last property is formulated as follows. Let
$(\ga,\a)$ be a $C^{*}$--dynamical system consisting of the
$C^{*}$--algebra $\ga$ and the automorphism $\a$. The unique
ergodicity for $(\ga,\a)$ is equivalent (cf. \cite{AD, AM, MT}) to
the norm convergence
\begin{equation}
\label{er}
\lim_{n\to+\infty}\frac1n\sum_{k=0}^{n-1}\a^{n}(a)=E(a)\,,\quad
a\in\ga\,.
\end{equation}
Here, $E=\om(\,\cdot\,)\idd$ is the conditional expectation onto the
fixed point subalgebra of $\a$ consisting of the constant multiples
of the identity, and $\om\in\cs(\ga)$ is the unique invariant state
for $\a$, here $S(\ga)$ denotes the set of all states of $\ga$. A
natural generalization requires that the the fixed point subalgebra
for $\a$ in \eqref{er} is nontrivial. This property, denoted as the
unique ergodicity w.r.t. the fixed point subalgebra, has been
investigated in \cite{AD,AM}. The unique weak mixing was
investigated in \cite{FM, MT}, which means that
\begin{equation*}
\lim_{n\to+\infty}\frac{1}{n}\sum_{k=0}^{n-1}\big|\f(\a^{k}(a))-\f(E(a))\big|=0\,,
\quad a\in\ga\,,
\end{equation*}
for every $\f\in\cs(\ga)$. As before, $E$ is the unique conditional
expectation projecting onto the fixed point subalgebra. Finally, the
unique mixing was defined and investigated in \cite{F}. We simply
require that
\begin{equation}
\label{mta2} \lim_{n}\f(\a^{n}(a))=\f(E(a))\,, \quad a\in\ga\,,
\end{equation}
for every $\f\in\cs(\ga)$.

The property \eqref{mta2} of the convergence to the equilibrium is
perfectly meaningful in the quantum setting but its classical
counterpart is the following: if a classical system fulfils
\eqref{mta2} with $E(f)=\int f\di\m$, the support of the unique
invariant measure $\m$ is a singleton, that is, it is conjugate to
the trivial one point dynamical system, see \cite{F}. Indeed, we can
exhibit some interesting examples of uniquely mixing
$C^*$--dynamical system in the quantum setting, for which the fixed
point algebra is trivial or non trivial as well. Such examples are
constructed by quantizing the shift on $\ell^2(\bz)$ on the
Boltzmann--Fock space $\cf(\ell^2(\bz))$. Other examples arise from
some generalizations of the free shift to the shift on the reduced
$C^{*}$--algebras of $RD$--groups, and the free amalgamated product
of $C^{*}$--algebras (cf. \cite{F, FM}), respectively. Among these
case, there are examples of dynamical systems which are uniquely
mixing w.r.t. the fixed point algebra (i.e. for which there are many
invariant states, see below) as well. Finally, in \cite{DF} it has
been established that the shift automorphism of the $q$--deformed
Canonical Commutation Relations algebra enjoys unique mixing
property. We can then exhibit uniquely mixing $C^*$--dynamical
system for which the von Neumman algebra generated by the GNS
representation of the unique invariant state is a type $I_{\infty}$
factor (case of the unital algebra $\gar$ acting on the
Boltzmann--Fock space $\cf(\ch)$ generated by the annihilators
$a(f)$, $f\in\ch$), or a type $II_1$ factor (case of the unital
algebra $\gg$ acting on the Boltzmann--Fock space generated by the
selfadjoint part of the annihilators $a(f)+a^+(f)$, $f\in\ch$). The
reader is referred to the papers \cite{AD, AM, DF, F, FM, M} for
further details on the topic. It is then natural to address the
possibility to exhibit further examples of $C^*$--dynamical systems
enjoying all the ergodic properties such that the von Neumann
algebra generated by the GNS representation of the unique invariant
state is a factor of different type from the previous ones. The aim
of the present paper is to show that this is indeed possible by
quantizing classical dynamical systems satisfying the corresponding
ergodic properties.

Let $\ch$ be a separable Hilbert space. Consider for $-1\leq q\leq1$
the $q$--canonical commutation relations for annihilators $a(f)$ and
creators $a^+(f)$:
\begin{equation}
\label{krel} a(f)a^+(g)-qa^+(g)a(f)=\langle
g,f\rangle_{\ch}\idd\,,\qquad f,g\in\ch\,.
\end{equation}
The case $q=-1$ is the Fermionic (Canonical Anticommutation
Relation) case, whereas the case $q=1$ is the Bosonic (Canonical
Commutation Relation) case, and finally $q=0$ corresponds to the
Boltzmann (or free) case. For unitaries $U$ acting on $\ch$, the
Bogoliubov automorphisms, defined as $\a_U(a(f)):=a(Uf)$ are widely
investigated (e.g. \cite{BR}) in the CAR and CCR cases for the
natural physical applications. The reader is referred to \cite{BG,
BDS, D, NS, SV}, for various results and applications, including the
computation of various kind of entropies of the Bogoliubov
automorphisms. In \cite{BC} (see also \cite{V1}) it has been pointed
out by a few words, that the free shifts of the Cuntz algebra
$\co_\infty$ and on $C^*_r(\bbf_{\infty})$ are "highly ergodic". In
the present paper (see also the previous one \cite{F}) the meaning
of the previuos sentence is then clarified: the Bogoliubov
automorphisms, including the shift, are highly ergodic automorphisms
as they enjoy one of the strong ergodic properties of unique
ergodicity, unique weak mixing or unique mixing described below.

We start from a classical dynamical system $(X, T, \m)$ consisting
of a probability space $(X, \m)$, and a measure preserving
invertible transformation $T:X\mapsto X$. By using the Shlyakhtenko
construction on the Boltzmann--Fock space (cf. \cite{S, H}), we
consider the Bogoliubov automorphism $\a_U$ relative to the unitary
$U$ associated to the measure preserving transformation $T$. We
obtain the following results. If the classical dynamical system
$(X,T, \m)$ is ergodic but not weakly mixing, then the resulting
quantized system $(\gg,\a)$ is uniquely ergodic w.r.t. the fixed
point algebra (which is always nontrivial in this situation), but
not uniquely weak mixing. If we quantize a classical system
$(X,T,\m)$ which is weakly mixing but not mixing, the resulting
quantized system is uniquely weak mixing but not uniquely mixing.
Finally, if we quantize a mixing system $(X,T, \m)$, the resulting
quantum system will be uniquely mixing. In such a way, it is
possible to exhibit uniquely weak mixing and uniquely mixing
$C^*$--dynamical systems whose von Neumann algebra generated by the
GNS representation associated to the unique invariant state is a
factor of type $I_{\infty}$, $II_{1}$ or $III_{\l}$ where
$\l\in(0,1]$.

For the sake of completeness, our results are extended to the
$C^*$--dynamical systems based on the $q$--commutation relations,
provided $|q|<\sqrt2-1$.

To end the present introduction, we point out few things. At the
classical level we have a wide class of uniquely ergodic dynamical
systems. The reader is referred to \cite{Da} and the literature
cited therein. In addition, starting from a measure preserving
ergodic dynamical system $(X,T,\m)$, it is possible to construct in
a canonical way a uniquely ergodic classical dynamical system
$(C(Y), \a_S)$ such that $(Y,S,\n)$ is conjugate to $(X,T,\m)$, $\n$
being the unique invariant probability measure on $Y$ invariant
under $S$.\footnote{The measure preserving transformations
$(X_j,T_j,\m_j)$, $j=1,2$, are said to be conjugate if there exist
$\m_j$--measurable sets $A_j\in X_j$ of full measure such that
$T_j(A_j)=A_j$, and a one--to--one measure preserving map
$S:A_1\mapsto A_2$ such that $T_2=S\circ T_1\circ S^{-1}$.} This is
nothing but the Jewett--Krieger Theorem (cf. \cite{Je, K}). For the
intermediate weak mixing situation nothing is yet known and the
Jewett--Krieger Theorem is not yet available in this case. Finally,
for the mixing case, any classical dynamical system enjoying
\eqref{mta2} with $E=\om(\,\cdot\,)\idd$ is conjugate to the one
point trivial dynamical system and then the Jewett--Krieger Theorem
cannot be carry out. Notice that our approach is similar to the
Jewett and Krieger one at least in principle. Namely, starting from
a classical dynamical system based on an measure preserving
transformation which is ergodic, weakly mixing or mixing, we can
construct in a functorial way, nontrivial quantum dynamical systems
(one for each type $I_{\infty}$, $II_{1}$ or $III_{\l}$,
$\l\in(0,1]$ of von Neumann factor) enjoying the unique ergodicity,
unique weak mixing or unique mixing, respectively.

\section{preliminaries}

In this section we recall some preliminaries concerning
$C^*$-dynamical systems.

Let $\ga$ be a $C^*$--algebra with unit $\idd$. By $\cs(\ga)$ we
denote the set of all states on $\ga$. For a (discrete)
$C^*$--dynamical system we mean a triplet $\big(\ga,\a,\om\big)$
consisting of a unital $C^*$-algebra $\ga$, an automorphism $\a$ of
$\ga$, and a state $\om\in\cs(\ga)$ invariant under the action of
$\a$. The pair $(\ga,\a)$ consisting of a unital $C^*$-algebra and
an automorphism as before is called a $C^*$--dynamical system as
well. Suppose now that the classical dynamical system $(X, T, \m)$
is merely based of a probability space $(X, \m)$, and a measure
preserving invertible transformation $T:X\mapsto X$. It is well
known that $T$ induces a unitary transformation acting on
$L^2(X,\m)$. Consider the natural restriction of $U$ to
$\ch:=L^2(X,\m)\ominus\bc1$, being $1$ the constant function. It is
known (see e.g. \cite{KSF}) that the dynamical system $(X, T, \m)$
is {\it ergodic, weak mixing, mixing} iff
\begin{equation}
\label{uter} \lim_{n\to\infty}\frac1n\sum_{k=1}^n\langle
U^k\xi,\eta\rangle=0\,,\quad\xi,\eta\in\ch\,,
\end{equation}
\begin{equation}
\label{utwm} \lim_{n\to\infty}\frac1n\sum_{k=1}^n|\langle
U^k\xi,\eta\rangle|=0\,,\quad\xi,\eta\in\ch\,,
\end{equation}
\begin{equation}
\label{utm} \lim_{n\to\infty}\langle
U^n\xi,\eta\rangle=0\,,\quad\xi,\eta\in\ch\,,
\end{equation}
respectively. In the present paper, we call with an abuse of
notations, any unitary $U$ acting on a Hilbert space $\ch$
satisfying \eqref{uter}, \eqref{utwm}, \eqref{utm}, {\it ergodic},
{\it weakly mixing}, or {\it mixing}, respectively.

Let $(\ga,\a)$ be a $C^*$--dynamical system, and $E:\ga\mapsto\ga$ a
linear map. suppose that
\begin{equation}
\label{mp1}
\lim_{n\to+\infty}\frac{1}{n}\sum_{k=0}^{n-1}\f(\a^{k}(x))=\f(E(x))\,,\quad
x\in\ga\,,\f\in\cs(\ga)\,,
\end{equation}
\begin{equation}
\label{mpp1}
\lim_{n\to+\infty}\frac{1}{n}\sum_{k=0}^{n-1}\big|\f(\a^{k}(x))-\f(E(x))\big|=0\,,\quad
x\in\ga\,,\f\in\cs(\ga)\,,
\end{equation}
or finally,
\begin{equation}
\label{mppp1} \lim_{n\to+\infty}\f(\a^{n}(x))=\f(E(x))\,,\quad
x\in\ga\,,\f\in\cs(\ga)\,.
\end{equation}

It can readily seen (cf. \cite{FM}) that the map $E$ is a
conditional expectation projecting onto the fixed point subalgebra
$\ga^{\a}:=\{x\in\ga: \ \a(x)=x\}$. Furthermore, $E$ is invariant
w.r.t. $\a$.
\begin{defin}
\label{smx} $(\ga,\a)$ is said to be uniquely ergodic, uniquely weak
mixing or uniquely mixing w.r.t the fixed point subalgebra, if
\eqref{mp1}, \eqref{mpp1} or \eqref{mppp1} holds true, respectively.

If $E=\om(\,\cdot\,)\idd$, then we simply call the dynamical system
$(\ga,\a)$ uniquely ergodic, uniquely weak mixing or uniquely mixing
(UE, UWM and UM for short), respectively.\footnote{If
$E=\om(\,\cdot\,)\idd$, then there is a unique invariant state for
$\a$, see \cite{AD}.}
\end{defin}
By using the Jordan decomposition of bounded linear functionals, one
can replace $\cs(\ga)$ with $\ga^*$ everywhere in Definition
\ref{smx}. We refer to \cite{St, SZ} for standard results on the
operator algebras and modular theory.

Let $\ch$ be a separable Hilbert space. The Boltzmann--Fock space
(called sometimes full Fock space) $\cf(\ch)$ is defined by
$$
\cf(\ch):=\bc\Om\oplus_{n=1}^{\infty}\ch^{\otimes n}\,.
$$
The vector $\Om$ is called {\it the vacuum vector}, and the vector
state $\om:=\langle\,\cdot\,\Om\,,\Om\rangle$ the vacuum state For
$f\in\ch$, the (left) creator $a^+(f)$ acts on $\cf(\ch)$ by
\begin{equation}
\label{qcre} 
a^+(f)\Omega=f\,,\quad a^+(f)f_1\otimes\cdots\otimes f_n=f\otimes
f_1\otimes\cdots\otimes f_n\,,
\end{equation}
and its adjoint is the (left) annihilator $a(f)$ given by
$$
a(f)\Omega=0\,,\quad a(f)f_1\otimes\cdots\otimes
f_n=\langle f_1,f\rangle f_2\otimes\cdots\otimes f_n\,.
$$
It is easily seen that $a(f)^*=a^+(f)$, and the $a(f)$ satisfy the
commutation rule \eqref{krel} with $q=0$. Let $P:\ch\mapsto\ch$ be a
positive contraction. Then the map $\cf(P):\cf(\ch)\mapsto\cf(\ch)$
defined as
$$
\cf(P)\Om:=\Om\,,\quad \cf(P)f_1\otimes\cdots\otimes
f_n:=Pf_1\otimes\cdots\otimes Pf_n
$$
is a contraction. The induced map $T_P(a(f))=a(Pf)$ is a completely
positive map on the unital $C^*$--algebra generated by the all the
annihilators. If $U$ is an isometry (resp. unitary), then $T_U$ is
an endomorphism (resp. automorphism). In the sequel we are
interested in the case when $U$ is unitary. Then $T_U$ is called a
{\it Bogoliubov automorphism} (see \cite{H, S, V}). Note that
dynamical and topological entropies of such kind of automorphisms
were intensively studied in \cite{BG, GN, SV, V1} for the CAR and
the CCR cases. The unique case relative to $q=0$ concerns the
quantization of the shift on $\ell^2(\bz)$, see e.g \cite{BC, BDS}.
Concerning the entropy, nothing is yet known for a general
Bogoliubov automorphism in the case $-1<q<1$.

\section{ergodic properties of Bogoliubov automorphisms}

In the present paper we assume that all the Hilbert spaces we deal
with are separable even if it is not directly specified.

Let $\ch$ be a separable Hilbert space and $U:\ch\to\ch$ be an
unitary. By $\a_U$ we denote the Bogoliubov automorphism $T_U$. We
start with the the following estimation needed in the sequel.
\begin{prop}
\label{mmain} Let $\{k_l\}_{l\in\bn}$ be any subsequence of natural
numbers, and $U$ be a unitary operator. Under the above notations,
we have the following estimation
\begin{align*}
&\bigg\|\sum_{l=1}^N\a_U^{k_l}\big(a^+(f_1)\cdots
a^+(f_m)a(g_1)\cdots
a(g_n)\big)\bigg\|\\
\leq&\bigg\|\sum_{l=1}^NU^{-k_l}f_1\otimes\cdots\otimes U^{-k_l}f_m
\otimes U^{k_l}g_n\otimes\cdots\otimes U^{k_l}g_1\big)\bigg\|\,.
\end{align*}
\end{prop}
\begin{pf}
It is enough to consider $x\in\ch^{\otimes t}$ with $t\geq n$. By
using any orthonormal basis $\{e_j\}_{j\in J}$ for $\ch$, we can
symbolically write
$$
x=\sum_{\s_1,\dots,\s_n,{\bf s}}x_{\s_1,\dots,\s_n,{\bf
s}}e_{\s_1}\otimes\cdots\otimes e_{\s_n}\otimes \xi_{{\bf s}}\,,
$$
with $\langle\xi_{{\bf r}},\xi_{{\bf s}}\rangle=\d_{{\bf r}{\bf
s}}$. We can also suppose that the $x_{\s_1,\dots,\s_n,{\bf s}}$ are
zero but finitely many of them. Put $F:=f_1\otimes\cdots\otimes
f_m$, $G:=g_n\otimes\cdots\otimes g_1$. We have
\begin{align*}
&\G:=\sum_{l=1}^N\a_U^{k_l}\big(a^+(f_1)\cdots a^+(f_m)a(g_1)\cdots
a(g_n)\big)x\\
=\sum_{l=1}^N\sum_{{\bf s}}\bigg\langle&\sum_{\s_1,\dots,\s_n}
x_{\s_1,\dots,\s_n,{\bf s}}e_{\s_1}\otimes\cdots\otimes
e_{\s_n},\big(U^{\otimes n}\big)^{k_l} G\bigg\rangle \big(U^{\otimes
m}\big)^{k_l}F\otimes\xi_{{\bf s}}\,.
\end{align*}
It follows that $\langle\G,\G\rangle$ can be viewed as a linear
combinations of inner products in $\ch^{\otimes(2n+m)}$, obtaining
\begin{align*}
&\langle\G,\G\rangle=\sum_{{\bf
s}}\bigg\langle\bigg(\sum_{\s_1,\dots,\s_n}
x_{\s_1,\dots,\s_n,{\bf s}}e_{\s_1}\otimes\cdots\otimes e_{\s_n}\bigg)\\
\otimes&\bigg(\sum_{l=1}^N\big(U^{\otimes
m}\big)^{-k_l}F\otimes\big(U^{\otimes n}\big)^{k_l}G \bigg),
\bigg(\sum_{l=1}^N\big(U^{\otimes
n}\big)^{k_l}G\otimes\big(U^{\otimes
m}\big)^{-k_l}F \bigg)\\
\otimes&\bigg(\sum_{\s_1,\dots,\s_n} x_{\s_1,\dots,\s_n,{\bf
s}}e_{\s_1}\otimes\cdots\otimes
e_{\s_n}\bigg)\bigg\rangle\\
\leq&\bigg\|\sum_{l=1}^N\big(U^{\otimes
m}\big)^{-k_l}\otimes\big(U^{\otimes n}\big)^{k_l} F\otimes G
\bigg\|^2\sum_{{\bf s}}\bigg\|\sum_{\s_1,\dots,\s_n}
x_{\s_1,\dots,\s_n,{\bf s}}e_{\s_1}\otimes\cdots\otimes
e_{\s_n}\bigg\|^2\\
=&\bigg\|\sum_{l=1}^N\big(U^{\otimes
m}\big)^{-k_l}\otimes\big(U^{\otimes n}\big)^{k_l} F\otimes G
\bigg\|^2\|x\|^2\,.
\end{align*}
\end{pf}

Let $\ch_\br$ be a separable real Hilbert space, and
$U_{\br}:\ch_{\br}\mapsto\ch_{\br}$ be an orthogonal transformation.
Extend $U_{\br}$ by linearity to $\ch_{\bc}:=\ch_{\br}+i\ch_{\br}$
and denote such a unitary operator as $U_{\bc}$. We report the
following known fact for the convenience of the reader.
\begin{lem}
\label{a} If $\s_{\mathop{\rm pp}}(U_{\bc})\neq\emptyset$, then
$U_{\br}\otimes U_{\br}$ has a nontrivial invariant vector.
\end{lem}
\begin{pf}
Let $e^{i\th}\in\s_{\mathop{\rm pp}}(U_{\bc})$ with eigenvector
$v=x+iy$, with $x,y\in\ch_{\br}$. Then
$$
U_{\br}\begin{pmatrix}
x\\
y\\
\end{pmatrix}=
\begin{pmatrix}
\cos\th &-\sin\th \\
\sin\th &\cos\th\\
\end{pmatrix}\begin{pmatrix}
x\\
y\\
\end{pmatrix}\,.
$$
The vector we are searching for is nothing but $x\otimes
x$+$y\otimes y$.\footnote{It follows that, if $e^{i\th}=\pm1$, then
$U_{\br}$ has an eigenvector corresponding to $e^{i\th}$. If
$e^{i\th}\in\bt\backslash\{\pm1\}$, then $U_{\br}$ has
$\begin{pmatrix}
\cos\th &-\sin\th \\
\sin\th &\cos\th\\
\end{pmatrix}$ as a direct summand.}
\end{pf}

Now consider a real subspace $\ck\subset\ch$ of the Hilbert space
$\ch$ and suppose that the unitary operator $U$ acting on $\ch$
satisfies $U\ck\subset\ck$, $U^*\ck\subset\ck$.\footnote{In this
case, $U\lceil_{\ck}$ defines an orthogonal transformation on $\ck$,
when the last is equipped with the inner product
$(x,y):=\text{Re}\langle x,y\rangle$.} Let $(\gar_{\ck},\a_U)$ be
the $C^*$--dynamical system, where $\gar_{\ck}$ is the unital
$C^*$--algebra acting on $\cf(\ch)$ generated by $\{a(f)\mid
f\in\ck\}$, and $\a_U$ the restriction of $T_U$ to $\gar_{\ck}$. The
$C^*$--dynamical system $(\gg_{\ck},\a_U)$ consists of the unital
$C^{*}$--algebra acting on $\cf(\ch)$ generated by
$\{s(f):=a(f)+a(f)^{+}\mid f\in\ck\}$ and the restriction of $T_U$
to $\gg_{\ck}$. Note that this is a Voiculescu's $C^*$--Gaussian
functor, see \cite{V}. In the sequel we simply write $(\gar,\a)$ and
$(\gg,\a)$, respectively.
\begin{prop}
\label{mmain1} If $U$ on $\ck$ is ergodic, then the dynamical
systems $(\gar,\a)$, $(\gg,\a)$ are ergodic w.r.t. the fixed point
algebra.
\end{prop}
\begin{pf}
By a standard approximation argument, it is enough to consider the
case when
\begin{equation}
\label{opprr} A:=a^+(f_1)\cdots a^+(f_m)a(g_1)\cdots a(g_n)\,,
\end{equation}
where the $f_1,\dots, f_m,g_1\dots,g_n$ are either eigenvectors of
$U$, or belong to the Hilbert subspace relative to the continuous
spectrum of $U$. If $f_1,\dots,f_m,g_1\dots,g_n$ are eigenvectors of
$U$ with corresponding eigenvalues
$e^{i\th_1},\dots,e^{i\th_m},e^{i\Th_1}\dots,e^{i\Th_n}$ such that
$$
\sum_{i=1}^m\th_i-\sum_{i=1}^n\Th_i=2h\pi
$$
for some integer $h\in\bz$, we conclude that $A$ is invariant under
$\a$. Otherwise, the vector $f_1\otimes\cdots\otimes f_m\otimes
g_1\otimes\cdots\otimes g_n\in\ch^{\otimes (m+n)}$ is not invariant
for the unitary $(U^*)^{\otimes m}\otimes U^{\otimes n}$. By
Proposition \ref{mmain} and the Mean Ergodic Theorem, we get in the
last case $\lim_N\frac1N\sum_{k=1}^N\a^k(A)=0$, and the proof
follows.
\end{pf}

\begin{prop}
\label{mmain2} If $U$ on $\ck$ is weakly mixing (resp mixing), then
the dynamical systems $(\gar,\a)$, $(\gg,\a)$ are UWM (resp. UM)
with the vacuum state the unique invariant state under $\a$.
\end{prop}
\begin{pf}
Let $A$ be as in \eqref{opprr}. We have for any subsequence
$\{k_l\}_{l\in\bn}$ of natural numbers of positive lower density,
$\lim_N\frac1N\sum_{l=1}^N\a^{k_l}(A)=0$ by taking into account
Proposition \ref{mmain} and the fact that $U$ (and then
$(U^*)^{\otimes m}\otimes U^{\otimes n}$) is weakly mixing (cf.
\cite{JL}). Again by \cite{JL}, this implies that, for each
$X\in\gar$ such that $\om(X)=0$, the sequence
$\{\a^n(X)\}_{n\in\bn}$ is (uniformly) weakly mixing at $0$. It
turns out to be equivalent to the fact that $(\gar,\a)$ is UWM with
$\om$ the unique invariant state. In the mixing case, by Proposition
\ref{mmain}, we have for operators $A$ as before and any subsequence
$\{k_l\}_{l\in\bn}$ of natural numbers,
$\lim_N\frac1N\sum_{l=1}^N\a^{k_l}(A)=0$ (cf. \cite{L}). The proof
follows by Proposition 2.3 of \cite{F}.
\end{pf}

Now we show that the quantized systems arising from ergodic but not
weakly mixing classical dynamical systems, cannot be UWM w.r.t. the
fixed algebra. The same will happen in the weak mixing situation:
the resulting quantum system cannot be UM.
\begin{prop}
\label{mmain3} Let $U$ be ergodic (resp weakly mixing) and suppose
that there exists some $f\in\ck$ such that the sequence
$\{U^kf\}_{k\in\bn}$ is not weakly mixing (resp. mixing) at $0$.
Then the dynamical systems $(\gar,\a)$, $(\gg,\a)$ cannot be UWM
(resp UM) w.r.t. the fixed point algebra.
\end{prop}

\begin{pf}
Let $f\in\ck$ such that $\{U^kf\}_{k\in\bn}$ is not weakly mixing
(resp. mixing) at $0$. By Proposition \ref{mmain1}, the dynamical
systems $(\gar,\a)$, $(\gg,\a)$ are UE w.r.t. the fixed point
algebra. Thus,
$$
\lim_N\frac1N\sum_{k=1}^N\a^k(a^+(f))=E(a^+(f))=0=E\lceil_\gg(s(f))\,,
$$
being $E$ the conditional projection onto $\gar^\a$. According to
\cite{JL} (\cite{L}) there exists a subsequence $\{k_l\}_{l\in\bn}$
of natural numbers of positive lower density (resp. a subsequence of
natural numbers) such that
$$
\limsup_N\bigg\|\frac1N\sum_{l=1}^NU^{k_l}f\bigg\|>0\,.
$$

Suppose that $(\gar,\a)$ or $(\gg,\a)$ is UWM (resp. UM) w.r.t. the
fixed point algebra. Then one gets
\begin{align*}
0=\lim_N\frac1N\bigg\|\sum_{k=1}^N\a^k&(s(f))\Om\bigg\|
=\lim_N\bigg\|\frac1N\sum_{k=1}^N\a^k(a^+(f))\Om\bigg\|\\
=&\limsup_N\bigg\|\frac1N\sum_{l=1}^NU^{k_l}f\bigg\|>0
\end{align*}
which is a contradiction.
\end{pf}

\section{on the type of the factors generated by Bogoliubov automorphisms}

In the present section we construct $C^*$--dynamical systems
enjoining the strong ergodic properties listed in Section 2, and
whose GNS representation relative to the Fock vacuum (which is the
unique invariant state for the discrete dynamics in the case of UWM
and UM) generates type $II_{1}$ and type $III_{\l}$, $\l\in(0,1]$
von Neumann factors. This is done by quantizing any classical
ergodic, weakly mixing or mixing dynamical system on the
Boltzmann--Fock space.

Let $(X, T,\m)$ be a classical dynamical system made of a
probability space $(X,\m)$, and a measure preserving invertible
transformation $T:X\mapsto X$. We suppose that $L^2(X,\m)$ is
infinite dimensional. Let
\begin{equation}
\label{eeqq}
\ck_{\br}:=(L^2_{\br}(X,\m)\ominus\br1)\bigotimes\big(\oplus_{\l\in
G}\br^2\big)\,.
\end{equation}
Here, $1\in L^2_{\br}(X,\m)$ is the constant f unction which is
invariant under the action of $U$, and $G$ is any countable
multiplicative subgroup of $\br_+$. Let $uf:=f\circ T^{-1}$ and
$$
v(t):=\bigoplus_{\l\in G}\begin{pmatrix}
\cos(t\ln\l) &-\sin(t\ln\l) \\
\sin(t\ln\l) &\cos(t\ln\l)\\
\end{pmatrix}\,.
$$
Then $u\otimes I$ and $I\otimes v(t)$ are orthogonal transformations
acting on the real Hilbert space $\ck_{\br}$ satisfying $[u\otimes
I,I\otimes v(t)]=0$. Let $\ck_{\bc}$ be the complexification of
$\ck_{\br}$ together with the positive non singular generator $A$ of
the complexification of $I\otimes v(t)$ as $I\otimes v(t)=I\otimes
a^{it}=:A^{it}$. Let $\ch$ be the completion of $\ck_{\bc}$ with
respect the inner product induced by $A$
\begin{equation}
\label{ccom} \langle x,y\rangle:=(2A(I+A)^{-1}x,y)\,,
\end{equation}
where $(\,{\bf\cdot}\,,\,{\bf\cdot}\,)$ is the inner product of
$\ck_{\bc}$. Denote by $U$ and $V(t)$ the unitary extension of the
corresponding orthogonal operators to the whole $\ch$. Let
$\cf(\ch)$ be the full Fock space generated by $\ch$ together with
the Fock vacuum vector $\Om$, and $\gg$ the $C^*$--algebra acting on
$\cf(\ch)$, generated by $\{s(f):=a(f)+a^+(f)\,:\,f\in\ck\}$. Notice
that $\Om$ is cyclic for $\gg$ and $\gg'$ (cf. \cite{S}), that is
$\Om$ is a standard vector for $\gg''$. The dynamical system under
consideration is $(\gg,\a)$, where $\a$ is the automorphism on $\gg$
induced by $\a(s(f)):=s(Uf)$.

\begin{prop}
\label{b} If $(X, T,\m)$ is ergodic but not weakly mixing, then the
fixed point algebra $\gg^{\a}$ is nontrivial.
\end{prop}
\begin{pf}
As $u$ is nontrivial and not weakly mixing, $U$ has at least an
eigenvalue $\chi$ in $\bt\backslash\{1\}$. If $\chi=-1$ there is a
corresponding eigenvector $f\in\ck$. An invariant element under the
action of $\a$ is $s(f)^2$. If $\xi\in\bt\backslash\{\pm1\}$, with
the corresponding eigenvector $v=f+ig$, then by Lemma \ref{a}, an
invariant element is given by $s(f)^2+s(g)^2$.
\end{pf}

The main results of the present paper are summarized in the
following theorems.
\begin{thm}
\label{mainnn} Let $\a$ be the automorphism in $\gg$ induced by
$\a(s(f)):=s(Uf)$. Then the following assertions hold true.
\begin{itemize}
\item[(i)] If $(X, T,\m)$ is ergodic but not weakly mixing, then
$(\gg,\a)$ is UE w.r.t $\gg^\a$, which is always nontrivial.
\item[(ii)] If $(X, T,\m)$ is
weakly mixing but not mixing, then $(\gg,\a)$ is UWM but not UM,
with $\om$ as the unique invariant state.
\item[(iii)] If $(X, T,\m)$ is
mixing, then $(\gg,\a)$ is UM, with $\om$ the unique invariant
state.
\end{itemize}
\end{thm}
\begin{pf}
We start by noticing that if $U$ is ergodic (resp. weakly mixing or
mixing) in $\ck_{\bc}$, then its extension on $\ch$ is ergodic
(resp. weakly mixing or mixing) as well. This easily follows by
\eqref{ccom} as, for any subsequence $\{k_l\}_{l\in\bn}$ of natural
number, we get
\begin{equation}
\label{nnoo}
\bigg\|\frac1N\sum_{l=1}^NU^{k_l}f\bigg\|_{\ch}\leq\|A\|^{1/2}_{\cb(\ck_{\bc})}\bigg\|\frac1N
\sum_{k=1}^NU^kf\bigg\|_{\ck_{\bc}}\,,
\end{equation}
$A$ being the positive operator in \eqref{ccom}. Then by Proposition
\ref{mmain1} (resp. Proposition \ref{mmain2}), $(\gg,\a)$ is UE
w.r.t. the fixed point algebra (resp. UWM or UM with the Fock vacuum
$\om$ as the unique invariant state). On the other hand, if $U$ is
not weakly mixing the pure point spectrum of $U$ is nonvoid.
Therefore the fixed point algebra $\gg^{\a}$ is nontrivial by
Proposition \ref{b}. Let now $F$ be a nonnull function on
$L^2(X,\m)$ with $\int F\di\m=0$ such that
\begin{equation}
\label{reim} \limsup_N\bigg\|\frac1N\sum_{l=1}^NU^{k_l}F\bigg\|>0
\end{equation}
for some subsequence $\{k_l\}_{l\in\bn}$ of natural numbers of
positive lower density (resp. a subsequence of natural numbers).
Notice that if $G=G_1+iG_2$, then
$\int|G|^2\di\m=\int(G_1^2+G_2^2)\di\m$. This means that
\eqref{reim} should be fulfilled at least by one of the real or
imaginary part of $F$. Thus, we can suppose without loss of
generality, that $F$ itself is real. Choose then $f:=F\otimes\xi$
with $\xi\in\oplus_{\l\in G}\br^2$. Then by \eqref{ccom} we have
$$
\limsup_N\bigg\|\frac1N\sum_{l=1}^NU^{k_l}f\bigg\|_{\ch}>0
$$
for the same subsequence $\{k_l\}_{l\in\bn}$ of natural numbers of
positive lower density (resp. a subsequence of natural numbers) as
before. Therefore, we conclude by Proposition \ref{mmain3} that if
$(X, T,\m)$ is ergodic but not weakly mixing (resp. weakly mixing
but not mixing), $(\gg,\a)$ cannot be UWM w.r.t. the fixed point
algebra (resp UM).
\end{pf}
\begin{thm}
For the $C^*$--dynamical systems considered above, we have that
$\gg''\cong\pi_{\om}(\gg)''$ is a non injective von Neumann factor
of type $II_1$, $III_{\lambda}$, $\lambda\in(0,1)$ or $III_1$,
whenever $G$ is $\{1\}$, $\{\l^n\,:\, n=0,1,2,\dots\}$ or $\bq_+$
respectively, $\pi_{\om}$ being the GNS representation relative to
$\om$.
\end{thm}
\begin{pf}
As we are assuming that $L^2(X,\m)$ is infinite dimensional, we have
that the positive operator $A$ in \eqref{ccom}, which is almost
periodic in our construction (cf. \cite{NSZ}), has always infinitely
many mutually orthogonal eigenvectors corresponding to the
eigenvalue $1$. Let $1=\l_1=\l_2=\cdots=\l_N=\cdots$ be an infinite
sequence of such eigenvalues. We have
$\frac1{\sqrt{N}}\sum_{k=1}^{N}\frac2{\sqrt{\l_k}+\sqrt{\l_k^{-1}}}>4$
whenever $N>16$. Therefore, $\gg''$ is not injective by Theorem 2.2
of \cite{H}. On the other hand, by Theorem 3.2 of \cite{H}, the
centralizer $(\gg'')_{\om}$ has trivial relative commutant in
$\gg''$. This implies that $\gg''$ is a factor. Finally, Theorem 3.3
of \cite{H} (see also \cite{S}) provides the result relative to the
type of the factor $\gg''$.
\end{pf}
Notice that $\gar''\cong\pi_{\om}(\gar)''$ is a type $I_{\infty}$
von Neumann factor, see e.g. \cite{DF} for the proof.

\section{the case of $q$--deformed commutation relations}

The present section is devoted to show that all the construction can
be carried out for the $q$--deformed commutation relations, at least
for sufficiently small $q$. For $-1<q<1$, the concrete
$C^*$--algebras $\gar_q$ and its subalgebra $\gg_q$ act on the
$q$--deformed Fock space $\cf_q(\ch)$, which is the completion of
the algebraic linear span of the vacuum vector $\Om$, together with
vectors
$$
f_1\otimes\cdots\otimes f_n\,,\quad
f_j\in\ch\,,j=1,\dots,n\,,n=1,2,\dots
$$
w.r.t. the $q$--deformed inner product
\begin{equation}
\label{qinn} 
\langle f_1\otimes\cdots\otimes
f_n\,,g_1\otimes\cdots\otimes g_m\rangle_q
:=\d_{n,m}\sum_{\pi\in\bp_n}q^{i(\pi)}\langle
f_1\,,g_{\pi(1)}\rangle\cdots\langle f_n\,,g_{\pi(n)}\rangle\,,
\end{equation}
$\bp_n$ being the symmetric group of $n$ elements, and $i(\pi)$ the
number of inversions of $\pi\in\bp_n$. The creator $a_q^+(f)$ is
defined as in \eqref{qcre}, and the corresponding annihilator is
defined as
\begin{align}
\label{qann}
a(f)\Omega&=0\,,\\
a(f)(f_1\otimes\cdots\otimes f_n)&=\sum_{k=1}^nq^{k-1}\langle
f_k,f\rangle f_1\otimes\cdots f_{k-1}\otimes
f_{k+1}\otimes\cdots\otimes f_n\,.\nn
\end{align}
$a_q^+(f)$ and $a_q(f)$ are adjoint each other w.r.t. the inner
product \eqref{qinn} and satisfy the commutation relations
\eqref{krel}. The Fock vacuum is defined as
$\om_q:=\langle\,{\bf\cdot}\,\Om,\Om\rangle$. As for the free
situation, a positive contraction $P$ on $\ch$ induces a completely
positive map $T_q(P)$ on the unital $C^*$--algebra generated by all
the annihilators $a_q(f)$ which is an isomorphism, provided $U$ is
unitary. The reader is referred to \cite{BKS} and the literature
cited therein, for further details.

Let $(\gg_q,\a_q)$ be the $C^*$--dynamical system where $\gg_q$ is
the unital $C^*$--algebra acting on $\cf_q(\ch)$ generated by
$\{s_q(f):=a_q(f)+a_q^+(f)\,|\,f\in\ck_{\br}\}$, and
$\a_q(s_q(f))=s_q(Uf)$. Here, $\ck_{\br}$ is given in \eqref{eeqq},
$\ch$ is the completion of $\ck_{\br}+i\ck_{\br}$ w.r.t. the inner
product given in \eqref{ccom}, and finally $U$ is the unitary acting
on $\ch$ as described in Section 4. 

In order to extend our previous results to the $q$--commutation relations, we need some preparatory results. 
\begin{prop}
\label{contt}
Let $\{U_n\}_{n\in\bn}$, $U$ be unitaries acting on $\ch$, together with the corresponding Bogoliubov automorphisms $\{\a^{(n)}_q\}_{n\in\bn}$ $\a_q$ on $\gar_q$, respectively. 

If $\lim_n U_n=U$ in the strong operator topology of $\cb(\ch)$, then $a^{(n)}_q$ converges pointwise in norm to $a_q$
\end{prop}
\begin{pf} 
It is enough to prove the assertion for each $a_q(f)$, $f\in\ch$. We get by Remark 1.2 of \cite{BKS},
\begin{align*}
\big\|\a^{(n)}\big(a_q(f)\big)-\a^{(n)}\big(a_q(f)\big)\big\|_{\cb(\cf_q(\ch))}&
=\|a_q(U_nf)-a_q(Uf)\|_{\cb(\cf_q(\ch))}\\
\leq&(1/\sqrt{1-|q|})\|(U_nf-Uf)\|_{\ch}
\end{align*}
\end{pf}

Let $\car$ be any finite dimensional Hilbert space whose dimension is equal to $d$, together with an orthonormal basis ${\displaystyle\{e_j\}_{j=1}^d}$. It is shown in 
\cite{JSW} and \cite{DN} that, if 
$|q|<\sqrt2-1$, the unital $C^*$--algebra $\gar_q$ is isomorphic to $\gar_0\equiv\gar$ via a map $\th$ sending $a_q(e_j)$ to $a_0(e_j)R$. Here, $R\in\gar_0$ is a positive element satisfying
$$
R^2=\sum_{j=1}^d a^+_0(e_j)a_0(e_j)+\sum_{j,k=1}^d (a_0(e_j)R 
a_0(e_k))^*(a_0(e_k)Ra_0(e_j))\,.
$$

Let 
${\displaystyle M:=\sum_{j=1}^d a^+_q(e_j)a_q(e_j)}$. It is a positive operator. Furthermore, 
consider the unitary operator 
$$ 
V=\bigoplus_{m=0}^{\infty}V_m:\cf_q(\car)\to\cf_0(\car)
$$ 
defined recursively as
$$
V_0:=J_0\,,\quad V_{m}:=(J_1\otimes V_{m-1})M^{1/2}\lceil_{\car^{\otimes m}}\,,\quad m=1,2,\dots\,\,,
$$
where $J_0$ $J_1$ are the identifications of the (complex multiple of the) vacuum vector and $\car$ in 
$\cf_q(\car)$, with the corresponding objects in $\cf_0(\car)$, respectively.
Then $\r$ can be written as $R=VM^{1/2}V^*$.
\begin{lem}
\label{finva}
Let $U$ be a unitary acting on the finite dimensional Hilbert space $\car$. Then $\cf_0(U)R=R\cf_0(U)$.
\end{lem}
\begin{pf}
By taking into account the definition of $M^{1/2}$, $V$ and $R$, it is enough to show that $M$ commutes with $\cf_q(U)$. We get
\begin{align*}
&\cf_q(U)Mf_1\otimes\cdots\otimes f_n\\
=&\sum_{i=1}^n q^{i-1}\bigg(\sum_{j=1}^d\langle f_i,e_j\rangle Ue_j\bigg)
\otimes Uf_1\otimes\cdots\otimes Uf_{i-1}\otimes Uf_{i+1}\otimes\cdots\otimes Uf_n\\
=&\sum_{i=1}^n q^{i-1}\bigg(\sum_{j=1}^d\langle Uf_i,Ue_j\rangle Ue_j\bigg)
\otimes Uf_1\otimes\cdots\otimes Uf_{i-1}\otimes Uf_{i+1}\otimes\cdots\otimes Uf_n\\
=&\sum_{i=1}^n q^{i-1}Uf_i
\otimes Uf_1\otimes\cdots\otimes Uf_{i-1}\otimes Uf_{i+1}\otimes\cdots\otimes Uf_n\\
=&M\cf_q(U)f_1\otimes\cdots\otimes f_n\,.
\end{align*}
\end{pf}

In addition, it is shown in \cite{H}, Section 5, that the previous result extends to the case of any separable Hilbert space $\car$, where $\th$ is the inductive limit of the corresponding isomorphisms 
$\th_n$ for each increasing sequence of $d_n$--dimensional subspaces $\car_n$ such that 
${\displaystyle\bigcup_n\car_n}$ is dense in $\car$. We refer the reader to the above mentioned paper \cite{H}, for further details relative to the isomorphism $\th$ realizing the equivalence between $\gar_q$ and $\gar_0$, when $|q|<\sqrt2-1$ (known in the literature as $E_q(\car)$ and 
$E_0(\car)$, respectively).
\begin{thm}
\label{iny}
There exists an isomorphism $\th:\gar_q\to\gar_0$  which  intertwines any Bogoliubov automorphism: 
$\th\circ\a_q=\a_0\circ\th$, provided $|q|<\sqrt2-1$.
\end{thm}
\begin{pf}
Let $U$ be the unitary acting on $\ch$ generating the Bogoliubov automorphism on the algebras 
$\gar_0$ and $\gar_q$, $|q|<\sqrt2-1$. Let $K$ be the Cayley transform of $U$, together with the resolution of the identity 
$\l\mapsto E(\l)$ of $K$, which is supposed to be right--continuous (in the strong operator topology).
Define the saw--tooth function
$$
h(\l):=\l-2k\pi\,,\quad\l\in(2k\pi,2(k+1)\pi]\,,\quad k\in\bz\,.
$$
 It is easy to show that $H:=\int\l\di \!E(\l)$ is a bounded selfadjoint operator such that $U=e^{iH}$. Fix an increasing sequence $\{\ch_n\}$ of finite dimensional subspaces such that 
 ${\displaystyle\bigcup_n\ch_n}$ is dense in $\ch$, togeheter with the associated selfadjoint projections $P_n$. Define $U_n:=e^{iP_nHP_n}$. We easily get $\lim_nU_n=U$ in the strong topology of $\car$. There exists an isomorphism 
 $\th:\gar_q\to\gar_0$ intertwining $\a^{(n)}_p$ and $\a^{(n)}_0$. This can be done by considering for each fixed $n$, the sequence of unitary operators $U_n\lceil_{\ch_m}$, $m>n$. Lemma 5.6 of \cite{H} leads to the claim by taking into account Lemma \ref{finva}. Fix $A\in\gar_q$. By Proposition 
 \ref{contt} we get
$$
\th(\a_q(A))=\th\big(\lim_n\a^{(n)}_q(A)\big)=\lim_n\th(\a^{(n)}_q(A))\lim_n\a^{(n)}_0(\th(A))
=\a_0(\th(A))\,,
$$
 which is the assertion.
\end{pf}
\begin{thm}
For the $C^*$--dynamical system $(\gg_q,\a_q)$, all the assertions
of Theorem \ref{mainnn} hold true, provided $|q|<\sqrt2-1$.
\end{thm}
\begin{pf}
Fix the $C^*$--dynamical system $(\gar_q,\a_q)$. By Theorem \ref{iny}. It is conjugate to the 
$C^*$--dynamical system $(\gar_0,\a_0)$, where $\a_q$, $\a_0$ are Bogoliubov automorphisms generated on $\gar_q$, $\gar_0$ by the same orthogonal operator. In addition, the isomorphism described in Proposition \ref{iny} intertwines the corresponding Fock vacua $\om_q$ and 
$\om_0$.\footnote{Notice that it is unclear if such isomorphism sends $\gg_q$ onto $\gg_0$. Thus, it is unclear if $(\gg_q,\a_q)$ is conjugate to $(\gg_0,\a_0)$.} By taking
into account of Proposition \ref{mmain} and \eqref{nnoo},
$(\gar_q,\a_q)$ is UE w.r.t. the fixed point subalgebra (resp UWM or
UM), provided $U$ is ergodic (resp. weakly mixing or mixing). Then
$(\gg_q,\a_q)$ is by restriction, UE w.r.t. the fixed point
subalgebra, UWM or UM, provided $U$ is ergodic, weakly mixing or
mixing, respectively. On the other hand, suppose that the classical
ststem $(X,T,\m)$ is ergodic (resp. weakly mixing) but not weakly
mixing (resp. mixing). We can choose as in the proof of Theorem
\ref{mainnn} $f\in\ck_{\br}$ such that
$$
\limsup_N\bigg\|\frac1N\sum_{l=1}^NU^{k_l}f\bigg\|_{\ch}>0
$$
for some subsequence $\{k_l\}_{l\in\bn}$ of natural numbers of
positive lower density (resp. a subsequence of natural numbers).
Then we have
\begin{align*}
0=\lim_N\frac1N\bigg\|\sum_{k=1}^N\a_q^k&(s_q(f))\Om\bigg\|
=\lim_N\bigg\|\frac1N\sum_{k=1}^N\a_q^k(a_q^+(f))\Om\bigg\|\\
=&\limsup_N\bigg\|\frac1N\sum_{l=1}^NU^{k_l}f\bigg\|>0
\end{align*}
which is a contradiction. Thus, $(\gg_q,\a_q)$ cannot be UWM w.r.t.
the fixed point subalgebra (resp. UM). Finally, in the case when
$(X,T,\m)$ is weakly mixing or mixing, the fixed point algebra of
$(\gar_q,\a_q)$ is trivial. Then by restriction, the fixed point
algebra of $(\gg_q,\a_q)$ is trivial as well. In the ergodic case,
the fixed point of $(\gg_q,\a_q)$ is non trivial by Proposition
\ref{b}.
\end{pf}

Consider now the $q$--shift. Namely, let
$u:\ell^2_{\br}(\bz)\mapsto\ell^2_{\br}(\bz)$ be the shift acting on
$\ell^2_{\br}(\bz)$, and
$$
\ck_{\br}:=\ell^2_{\br}(\bz)\bigotimes\big(\oplus_{\l\in
G}\br^2\big)\,.
$$
$\ch$ will be the completion of $\ck_{\br}+i\ck_{\br}$ w.r.t. the
inner product given in \eqref{ccom}. Define $U$ and $V(t)$ as in
Section 4, and finally $\gg_q$ as at the beginning of the present
section. In this case (as well as in all the cases previously
described in the present section for each $-1<q<1$) $\gg_q''$ is non
injective von Neumann factor of type $II_1$, $III_{\lambda}$,
$\lambda\in(0,1)$ or $III_1$, whenever $G$ is $\{1\}$, $\{\l^n\,:\,
n=0,1,2,\dots\}$ or $\bq_+$ respectively.

In \cite{DF}, it was proven for each $-1<q<1$, that the $q$--shift
is UM in the case when the modular theory is trivial (i.e. when
$G=\{1\}$ in \eqref{eeqq}). The same proof of Theorem 3 of \cite{DF}
allows us to extend the previous results to all cases $-1<q<1$, at
least in the case of the shift. We have then proven the following
\begin{prop}
For each $-1<q<1$, the $C^*$--dynamical system $(\gg_q,\a_q)$
($\a_q$ being the $q$--shift) is UM, with the Fock vacuum $\om_q$
the unique invariant state.
\end{prop}

We can conjecture that all the results described in the present
section for the $C^*$--dynamical systems based on the
$q$--commutation relations, can be extended to all $q\in(0,1)$.
Unfortunately, it is not known if all the $C^*$--algebras $\gar_q$
are isomorphic for any $q\in(-1,1)$. In addition, an estimation
similar to Proposition \ref{mmain} and to Proposition 2 of \cite{DF}
is not yet available for the general case when $q\neq0$ and the
involved Bogoliubov automorphism is not the shift.

\section*{Acknowledgement} The second--named author (F. M.) thanks the
MOHE grant FRGS0308-91.

\end{document}